\documentclass[12pt,a4paper]{article}
\usepackage{url}
\usepackage{graphicx}
\usepackage{amsmath}

\DeclareSymbolFont{AMSb}{U}{msb}{m}{n}
\DeclareMathSymbol{\N}{\mathbin}{AMSb}{"4E}
\DeclareMathSymbol{\Z}{\mathbin}{AMSb}{"5A}
\DeclareMathSymbol{\R}{\mathbin}{AMSb}{"52}
\DeclareMathSymbol{\Q}{\mathbin}{AMSb}{"51}
\DeclareMathSymbol{\I}{\mathbin}{AMSb}{"49}
\DeclareMathSymbol{\C}{\mathbin}{AMSb}{"43}

\newcommand{\Ci}{\operatorname{Ci}}
\newcommand{\Si}{\operatorname{Si}}
\newcommand{\si}{\operatorname{si}}
\newcommand{\res}{\operatorname{Res}}

\title{Single Exponential Approximation of Fourier Transforms}
\author{Patrick McLean\thanks{School of Mathematics and Physics, University of Tasmania, Private Bag 37, Hobart, Tasmania, 7001, \textsc{Australia}. \protect\url{mailto:p_mclean@maths.utas.edu.au}}}

\date{(31 November 2005)}

\begin{document}

\maketitle

\begin{abstract}
This article is concerned with a new method for the approximate evaluation of Fourier sine and cosine transforms. We develope and analyse a new quadrature rule for Fourier sine and cosine transforms involving transforming the integral to one over the entire real line and then using the trapezoidal rule in order to approximate the transformed integral. This method follows on from the work of Ooura and Mori, see \cite{om1} and \cite{om2}

A complete error analysis is made using contour integration. An example is examined in detail and the error is analysed using residues and the saddle point method. The method we have developed is characterised by its simplicity and single exponential convergence. 
\end{abstract}

\tableofcontents

\section{Introduction}
In this paper we consider the numerical approximation of Fourier sine and cosine transforms, that is, integrals of the form
\begin{eqnarray}
f_c(t)& = & \int_0^\infty f(x) \cos(tx) dx, \\
f_s(t)& = & \int_0^\infty f(x) \sin(tx) dx. \label{fst}
\end{eqnarray}
There exist extensive tables of Fourier transforms, see, for example, \cite{erdelyi}. Nonetheless there is need for the numerical approximation of Fourier transforms. We focus on \( f_s(t) \) and present a new method involving transforming the integral (\ref{fst}) to one over \( (-\infty,\infty) \) and then using the trapezoidal rule.

In section 2 we review the trapezoidal rule on \( (-\infty,\infty) \)and sources of error in its use. Quadrature rules over intervals other than \( (-\infty,\infty) \) arise from the use of transformations. The double exponential and sinc methods of quadrature arise in this way, see \cite{takahasimori} and \cite{stenger}, respectively.

In section 3 we introduce quadrature methods for Fourier sine integrals based on transformations of a certain sort. Double exponential versions were introduced in \cite{om1} and \cite{om2}. Here we present a single exponential version together with useful asymptotic estimates of its behaviour. We also note that the midpoint rule may be used for Fourier cosine transforms.

In section 4 we present an example and analyse the discretisation error using residues and the saddle point method. Here we see that the rate of convergence of our method is determined by the proximity of singularities of the integrand to the interval of integration \( (0,\infty) \). This phenomena is common in quadrature, see \cite{hough}, for example .

In section 5 we look at the truncation error and compare our method with the double exponential method of Ooura and Mori.

\section{The Trapezoidal Rule and its Error}
Given a function \( F(u) \) defined on \( (-\infty,\infty) \) with integral
\begin{equation}
I = \int_{-\infty}^\infty F(u) \, du,
\end{equation}
and a real number \( h>0 \) we define the trapezoidal approximation \( T_h \) to \( I \) by
\begin{equation}
\label{trap}
T_h = h \sum_{k=-\infty}^\infty F(kh),
\end{equation}
This rule is the basis of Sinc methods extensively developed by Stenger and can be obtained by integrating a Sinc function interpolant to \( F(u) \), see \cite{stenger}.

\subsection{Discretisation Error of the Trapezoidal Rule}
The quantity \( I-T_h \) is referred to as discretisation error. We shall use the error representation of Donaldson and Elliott \cite{de}. If we set \( \lambda = 0 \), \(a=0 \) and \( \nu = 1/h \) then equations (6.1) and (6.4) of \cite{de} read
\begin{eqnarray}
\label{trappsi}
\Psi_h(w) & = &\begin{cases} \pi \exp(i \frac{\pi}{h} w), & \Im w>0\\ \pi \exp(-i \frac{\pi}{h} w), & \Im w<0 \end{cases}\\
\label{trapphi}
\Phi_h(w) & = & - \sin(\frac{\pi}{h} w),
\end{eqnarray}
respectively.
Thus, from \cite[eqn 2.4]{de} we have the error representation
\begin{equation}
\label{traperr}
I-T_h = \frac{1}{2 \pi i} \int_C \frac{\Psi_h(w)}{\Phi_h(w)} F(w) \, dw,
\end{equation}
where \( C \) is a positively described contour enclosing the zeroes of \( \Phi_h(w) \)and avoiding any singularities of \( F(w) \).

\subsection{Truncation Error of the Trapezoidal Rule}
In practise, the (infinite) trapezoidal sum \( T_h \) is truncated at some \( \pm n \) to give the finite trapezoidal sum \(T_{n,h} \):
\begin{equation}
\label{trapfin}
T_{n,h} = h \sum_{k=-n}^n F(kh),
\end{equation}
as an approximation to \( I \). This introduces an error \( T_h-T_{n,h} \) that we refer to as truncation error. 

Thus, the overall error can be seen as comprising two sources:
\begin{equation}
I-T_{n,h}=I-T_h+T_h -T_{n,h}
\end{equation}
It is essential that the discretisation error and the truncation error match in order to achieve a suitable rate of convergence.

\section{Trapezoidal rule methods for Fourier sine transforms}

In this section we propose a method of approximation of Fourier integrals which is highly accurate and applicable to a wide range of functions including ones with singularities and slow decay. Our approximation to \( f_s(t) \) is obtained by making a particular change of variable from \( (0,\infty) \) to \( (-\infty,\infty) \), followed by an application of the trapezoidal rule. For \( f_c(t) \) we apply the same change of variable, but use the midpoint rule.

In order to approximate \( f_s(t) \) we take the integral ~(\ref{fst}) and introduce the change of variable 
\begin{equation}
x= \frac{m}{t} \phi(u), \qquad -\infty<u<\infty,
\end{equation}
where \( m>0 \) and \( \phi:(-\infty,\infty)\rightarrow (0,\infty) \) is a function satisfying 
\begin{eqnarray}
\label{as-}
\phi(u) \sim 0 \qquad \text{as} \qquad u\rightarrow -\infty, \\
\label{as+}
\phi(u) \sim u \qquad \text{as} \qquad u\rightarrow +\infty.
\end{eqnarray}
Thus, \( f_s(t) \) can be represented by the integral
\begin{equation}
\label{int1}
f_s(t) = \int_{-\infty}^\infty F_m(u) \, du
\end{equation}
where the function \( F_m(u) \) is defined by
\begin{equation}
\label{F_m}
F_m(u) = f(\frac{m}{t}\phi(u)) \sin(m\phi(u)) \frac{m}{t} \phi'(u).
\end{equation}
We now apply the trapezoidal rule with stepsize \( h=\frac{\pi}{h} \) to (\ref{int1}) to give the approximation
\begin{equation}
T_m=\frac{\pi}{m} \sum_{k=-\infty}^\infty F_m(k\frac{\pi}{m}).
\end{equation}
The asymptotic behaviour of \( \phi(u) \) at \( \pm \infty \) allows us to truncate this series without incurring too large an error. We shall provide further details on this in the last section of this paper.

\subsection{Midpoint rule for Fourier cosine transforms}
For completeness we note that \( f_c(t) \) may be efficiently evaluated by making the transformation \( x= \frac{m}{t} \phi(u) \) and then applying the midpoint rule\begin{equation}
M_h=h\sum_{k=-\infty}^\infty F((k+\frac 1 2)h),
\end{equation}
with stepsize \( h=\frac{\pi}{m} \).
The analogous error representation is 
\begin{equation}
\label{miderror}
I-M_h = \frac{1}{2 \pi i} \int_{C} \frac{\Psi_h(w)}{\Phi_h(w)} F_h(w)\, dw,
\end{equation}
where \( \Psi_h(w) \) and \( \Phi_h(w) \)  are given by 
\begin{eqnarray}
\label{midpsi}
\Psi_h(w) & = &\begin{cases} -i \pi \exp(i \frac{\pi}{h} w), & \Im w>0\\ i \pi \exp(-i \frac{\pi}{h} w), & \Im w<0 \end{cases}\\
\label{midphi}
\Phi_h(w) & = &\cos(\frac{\pi}{h} w),
\end{eqnarray}
respectively. Further details including an example are provided in \cite{thesis}.

\subsection{Double exponential transformation I}
In their first paper \cite[eq 11]{om1} Ooura and Mori introduce the transformation
\begin{equation}
\phi_1(u) = \begin{cases} \frac{u}{1-\exp(-K \sinh \, u)}, & u \neq 0\\ \frac{1}{K} & u=0\end{cases}.
\end{equation}
for \( K>0 \). They make the choice of \( K=2 \pi \).

\subsection{Double exponential transformation II}
In their second paper \cite[eqn 3.3]{om2} Ooura and Mori introduce the transformation
\begin{equation}
\phi_2(u) = \frac{u}{1-\exp(-2 u - \alpha(1-e^{-u})-\beta(e^u-1))}
\end{equation}
where the parameters are given by
\begin{equation}
\alpha= \beta/\sqrt{1+M \log(1+M)/(2 \pi)}, \quad \beta = \frac{1}{4}.
\end{equation}

\subsection{Novel Single Exponential Transformation}
We define the transformation \( \phi:(-\infty,\infty) \rightarrow (0,\infty) \) by
\begin{equation}
\label{trans}
\phi(u) = \log( e^u+1), \qquad -\infty < u < \infty,
\end{equation}
with derivative
\begin{equation}
\phi'(u) = \frac{e^u}{e^u+1}.
\end{equation}
Note that as a function of the complex variable \( w=u+iv \), \( \phi(w) \) has branch points at \( w= (2k+1)i \pi \) for an integer \( k \).

\subsubsection{Asymptotic behaviour of transformation}
We shall need to know the asymptotic behaviour of \( \phi(u) \) as \( u\rightarrow \pm \infty \).
The Taylor series of \( \log(1+x) \) around \( x= 0 \) is
\begin{equation}
\log(1+x) = x + O(x^2),\text{   for   } |x|<1,
\end{equation}
see \cite[eqn 4.1.24]{as}. Thus, we have that
\begin{equation}
\label{assminf}
\phi(u) = 0 + e^u +O(e^{2u}),
\end{equation}
as \( u\rightarrow -\infty \), and
\begin{equation}
\label{asspinf}
\phi(u) = u+\log(1+e^{-u}) = u + e^{-u} + O(e^{-2u}),
\end{equation}
as \( u \rightarrow \infty \).

\subsubsection{Asymptotic behaviour of inverse}
The inverse \( \phi^{-1}:(0,\infty)\rightarrow ( -\infty,\infty) \) is defined by
\begin{equation}
\label{transinv}
\phi^{-1}(x) = \log(e^x-1), \qquad 0 < x < \infty.
\end{equation}

We shall determine the asymptotic behaviour of \( w = \phi^{-1}(z) \) as \( z \rightarrow 0 \).
From (\ref{transinv}) we have
\begin{equation}
w= \log(e^z-1)  =\log z +\log\frac{e^z-1}{z}.
\end{equation}
On expanding the last term about \(z = 0 \), we arrive at
\begin{equation}
\label{phiinvat0}
\phi^{-1}(z) \sim \log z +\frac{z}{2}+O(z^2),
\end{equation}
as \( z \rightarrow 0 , \arg z \ne \pi \).

\section{A Certain Integral}
In this section we shall consider the integral
\begin{equation}
I = \int_0^\infty \frac{\sin tx}{(x-a)^2+b^2} \, dx,
\end{equation}
for \( t>0 \), \( -\infty<a<\infty \) and \( b>0 \). We shall estimate the discretisation error using residues and the saddle point method. We shall see that the discretisation error depends on the proximity of the singularities of the integrand \( a \pm i b \) to the interval of integration \( (0,\infty) \). 

\subsection{Evaluation in Terms of Trigonometric Integrals}
We follow the standard reference work \cite[\S 5.2]{as} and define the trigonometric integrals as follows.
For \( z \in \C \) define the sine integral \( \Si(z) \) as
\begin{equation}
\Si(z) = \int_0^z \frac{\sin t}{t} dt.
\end{equation}
The sine integral \( \Si \) is an entire function. A commonly used notation is \(\si(z) = \frac{\pi}{2} -\Si(z) \).

For \( z \in \C \) such that \( | \arg z | < \pi \) define the cosine integral \( \Ci(z) \) as
\begin{equation}
\Ci(z) = \gamma + \log z + \int_0^z \frac{\cos t-1}{t} dt.
\end{equation}
The cosine integral \( \Ci \) has a branch cut discontinuity along the negative real axis.

The functions \( \Si \) and \( \Ci \) occur as Fourier sine and cosine transforms.
Specifically, for \( |\arg a| < \pi \) and \( y>0 \) we have that
\begin{eqnarray}
\label{eq:cosint}
\int_0^\infty \frac{\cos(xy)}{a+x} \, dx & = & -\si(ay)\sin(ay) -\Ci(ay)\cos(ay) \\ 
\label{eq:sinint}
\int_0^\infty \frac{\sin(xy)}{a+x} \, dx & = & \Ci(a y)\sin(ay) -\si(a y)\cos(ay),
\end{eqnarray}
see \cite[\S 1.1 (9)]{erdelyi} and \cite[\S 1.2 (10)]{erdelyi}, respectively.

Using partial fractions together with ~(\ref{eq:cosint}) and ~(\ref{eq:sinint}) it is possible to show that
\begin{eqnarray}
bI &=& \sin at\sinh bt\Re\Ci(-a+ib)t+\cos at\cosh bt\Im\Ci(-a+ib)t\\
& & -\cos at\sinh bt \Re\si(-a+ib)t-\sin at\cosh bt\Im\Si(-a+ib)t).
\label{I}
\end{eqnarray}

\subsection{Discretisation error}
In order to implement our method we make the change of variable \( x=\frac{m}{t}\phi(u) \) and apply the trapezoidal rule with stepsize \( \frac{\pi}{m} \) to give the approximation
\begin{equation}
T_m = \frac{\pi}{m} \sum_{k=-\infty}^\infty F_m(k\frac{\pi}{m}),
\end{equation}
where the function \( F_m(u) \) is given by
\begin{equation}
F_m(u) = \frac{\sin(m\phi(u)) \frac{m}{t} \phi'(u)}{(\frac{m}{t}\phi(u)-a)^2+b^2}.
\end{equation}
We shall denote the pole of \( F_m(w) \) in the upperhalf plane closest to the real axis by \( w_0 \) and we note that \( \overline{w_0} \) is also a pole of \( F_m(w) \).

\begin{figure}[tbp]
\includegraphics{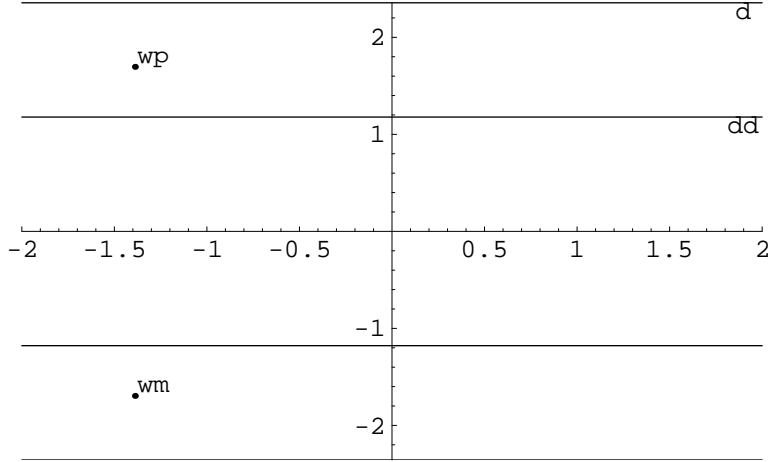}
\caption{The contours \( C \) and \( C' \)\label{poles}}
\end{figure}

Now, from (\ref{traperr}) we have that
\begin{equation}
\label{error}
I-T_m = \frac{1}{2 \pi i} \int_C \frac{\Psi_m(w)}{\Phi_m(w)} F_m(w)\, dw,
\end{equation}
where \( \Psi_m(w) \) and \( \Phi_m(w) \) are given by(\ref{trappsi}) and (\ref{trapphi}) respectively, and where \( C \) is a positively described closed contour going between the real axis and the poles \( w_0 \) and \( \overline{w_0} \) as depicted in figure \ref{poles}.

If we deform the contour to contain the poles \( w_0, \overline{w_0} \) we have by Cauchy's theorem
\begin{equation}
\label{error1}
I-T_m =  -\res(\frac{\Psi_m}{\Phi_m} F_m; w_0)-\res(\frac{\Psi_m}{\Phi_m} F_m;\overline{w_0}) + \frac{1}{2 \pi i} \int_{C'} \frac{\Psi_m(w)}{\Phi_m(w)} F_m(w)\, dw,
\end{equation}
where \( C' \) is a positively described contour going between the poles \( w_0, \overline{w_0} \) and the singular points of \( \phi(w) \) as depicted in figure \ref{poles}. We shall denote the two residue terms by \( R_m \) and the integral term by \( S_m \).

\subsection{Contribution to the error from Poles}
In this section we consider the evaluate the residue terms in ~(\ref{error1}).

Now, on using L'Hopital's rule, the residue of \( F_m(w) \) at \( w_0 \) is given by
\begin{eqnarray}
\res(F_m;w_0)&  =& \lim_{w\rightarrow w_0} (w-w_0) \frac{\sin(m\phi(w)) \frac{m}{t} \phi'(w)}{(\frac{m}{t}\phi(w)-a)^2+b^2} \\
& = & \frac{\sin(m \phi(w_0))} {2 (\frac{m}{t} \phi(w_0)-a) }.
\end{eqnarray}
Now, we use the fact that \( m \phi(w_0) = (a+ib)t \) to give
\begin{equation}
\res(F_m;w_0) = \frac{\sin((a+ib)t)}{2ib}.
\end{equation}

Similarly, we have that the residue of \( F_m(w) \) at \( \overline{w_0} \) is given by
\begin{equation}
\res(F_m;\overline{w_0}) = \frac{\sin((a-ib)t)}{-2ib}.
\end{equation}

Combining these residues with the definitions of \( \Psi_m(w) \) and \( \Phi_m(w) \) we have that
\begin{eqnarray}
R_m & = & -\frac{\pi \exp(i m w_0)}{-\sin(i mw_0)}  \frac{\sin(a+ib)t}{2ib} -\frac{\pi \exp(-i m\overline{w_0})}{-\sin(i m\overline{w_0})}  \frac{\sin(a+ib)t}{-2ib}\\
& = & \frac{\pi}{b} \left[ \frac{\sin(a+ib)t}{1-\exp(-2m i w_0)}+ \frac{\sin(a-ib)t} {1-\exp(2m i \overline{w_0})} \right]. 
\end{eqnarray}

We now substitute \( w_0=u_0+iv_0 \) and \( \overline{w_0} =u_0-iv_0 \) and rearrange to give
\begin{equation}
R_m = \frac{\pi}{b} \frac{ \left[ \exp(-2mv_0) - \cos(2mu_0)\right] \sin(at)\cosh(bt) + \sin(2mu_0)\cos(at)\sinh(bt)}{\cosh(2mv_0) - \cos(2mu_0)}.
\label{exact}
\end{equation}
Now, on neglecting the cosine term in the denominator and the exponential terms \( \exp\left(-2m v_0 \right) \) in the numerator we have that
\begin{equation}
R_m \sim \frac{2\pi}{b} \frac{ - \cos(2mu_0)\sin(at)\cosh(bt) +\sin(2mu_0) \cos(at) \sinh(bt)}{\exp(2mv_0)+\exp(-2mv_0)}
\end{equation}

On neglecting the term \( \exp \left(-2m v_0 \right) \) in the denominator, we obtain
\begin{equation}
R_m \sim \frac{2\pi}{b} \exp(-2m v_0 ) \left[ -\cos (2mu_0) \sin(at) \cosh(bt) +\sin (2mu_0)  \cos(at) \sinh(bt) \right].
\label{asympt}
\end{equation}

We shall only be interested in the asymptotic behaviour of \( w_0 \) as \( m \rightarrow \infty \), which we can obtain from (\ref{phiinvat0}):
\begin{equation}
\label{w0sim}
w_0 =\log\frac{(a+ib)t}{m} +\frac{(a+ib)t}{2m} +O(m^{-2}),
\end{equation}
as \( m \rightarrow \infty \). Now since \( \log(re^{i\theta}) = \log(r)+i \theta \) for \( r>0, 0<\theta <\pi \), it follows that the imaginary part \( v_0 \) of \( w_0 \) has asymptotic behaviour
\begin{equation}
\label{v0sim}
v_0 = \begin{cases} \arctan(\frac{b}{a}) +\frac{bt}{2m} +O(\frac{1}{m^2}), & a>0 ,\\
\frac{\pi}{2}+\frac{bt}{2m}, & a=0 \\
\pi-\arctan(\frac{b}{a}) +\frac{bt}{2m} +O(\frac{1}{m^2}), & a<0.
\end{cases}
\end{equation}
as \( m \rightarrow \infty \).

Thus using the asymptotic estimates (\ref{asympt}) and (\ref{v0sim}) we have the estimate
\begin{equation}
R_m \sim \begin{cases} C \exp \left(-2m \arctan(\frac{b}{a}) \right) , & a>0 ,\\
C \exp \left(-2m \frac{\pi}{2}) \right), & a=0 \\
C \exp \left(-2m (\pi-\arctan(\frac{b}{a})) \right), & a<0,
\end{cases}
\label{est1}
\end{equation}
as \( m \rightarrow \infty \), where \( C \) represents a number independent of \( m \).

\subsection{Contribution to the error from the saddle points}
In order to estimate the contribution to the error \( I-T_m \) from the integral \( S_m \) in (\ref{error1}) we write its integrand as \( \exp(p(w))q(w) \) where 
\begin{eqnarray}
p(w) &=& \log( \frac{\Psi_m(w)}{\Phi_m(w) }\sin(m\phi(w) )\\
q(w) & = & f(\frac{m}{t}\phi(w) ) \frac{m}{t} \phi'(w).
\end{eqnarray}
From the definitions of \( \Psi_m(w) \) and \( \Phi_m(w) \) we have that
\begin{equation}
p(w) = \log\left(-\pi(\cot(mw)+i)\sin(m\phi(w)) \right), \qquad \Im(w) >0.
\end{equation}
The saddle points of \( p(w) \) are solutions of the equation \( p'(w) =0 \). We shall be interested in the one in the upper half plane closest the real axis and we shall denote it by \( w_1 \). The point \( \overline{w_1} \) in the lower half plane is also a saddle point of \( p(w) \).

To find \( w_1 \) we must solve
\begin{equation}
p'(w) = -m(\cot(mw)-i)+\cot(m\phi(w))m \phi'(w)=0,
\end{equation}
that is, 
\begin{equation}
\label{problematic1}
\phi'(w)  \cot(m \phi(w)) = \cot(mw)-i.
\end{equation}
Solving this equation is problematic. Possible strategies are i) obtaining an asymptotic estimate for \(w_1 \) in \( m \) similar to the asymptotic estimate \( w_0 \) given by (\ref{w0sim}), or ii) given values of \( m \) numerically solving for \( w_1 \). Lacking the former we use the latter.

From \cite[eqn2.7.2]{debruijn} we have 
\begin{equation}
\frac{1}{2 \pi i} \int_{C_+} \frac{\Psi_m(w)}{\Phi_m(w)} F_m(w)\, dw
 \sim \frac{1}{2\pi i}\sqrt{\frac{2 \pi }{p''(w_1)}} \alpha e^{p(w_1)} q(w_1)
\label{saddlearse1}
\end{equation}
where \( C_+ \) is the upper half of \( C' \) and where 
\begin{equation}
\alpha= \exp(\frac{\pi i}{2}-\frac{i}2 \arg(p''(w_1))).
\end{equation}
The contribution to the error from the saddle point at \( \overline{w_1} \) is the conjugate of the estimate in (\ref{saddlearse1}). 
Hence, we have the estimate
\begin{equation}
S_m \sim 2\Re \left[ \frac{1}{2\pi i}\sqrt{\frac{2 \pi }{p''(w_1)}} \alpha e^{p(w_1)} q(w_1) \right],
\label{arse1}
\end{equation}
as an estimate of the integral over \( C' \) (see (\ref{error})).

\subsection{Discretisation error}

In this section we present results comparing the discretisation error \( I-T_m \) with the contribution from the poles \( R_m \) and the contribution from the saddle points \( S_m\). 

To evaluate \( I \) we use equation ~(\ref{I}) and to evaluate \( T_m \) to within machine precision we sum the series for sufficiently high values of \( n \). Typically taking \( n=4m^2 \). Furthermore, we use the representation  ~(\ref{exact}) for \( R_m \) and the estimate ~(\ref{arse1}) for \( S_m \). 

\begin{figure}[tbp]
\includegraphics{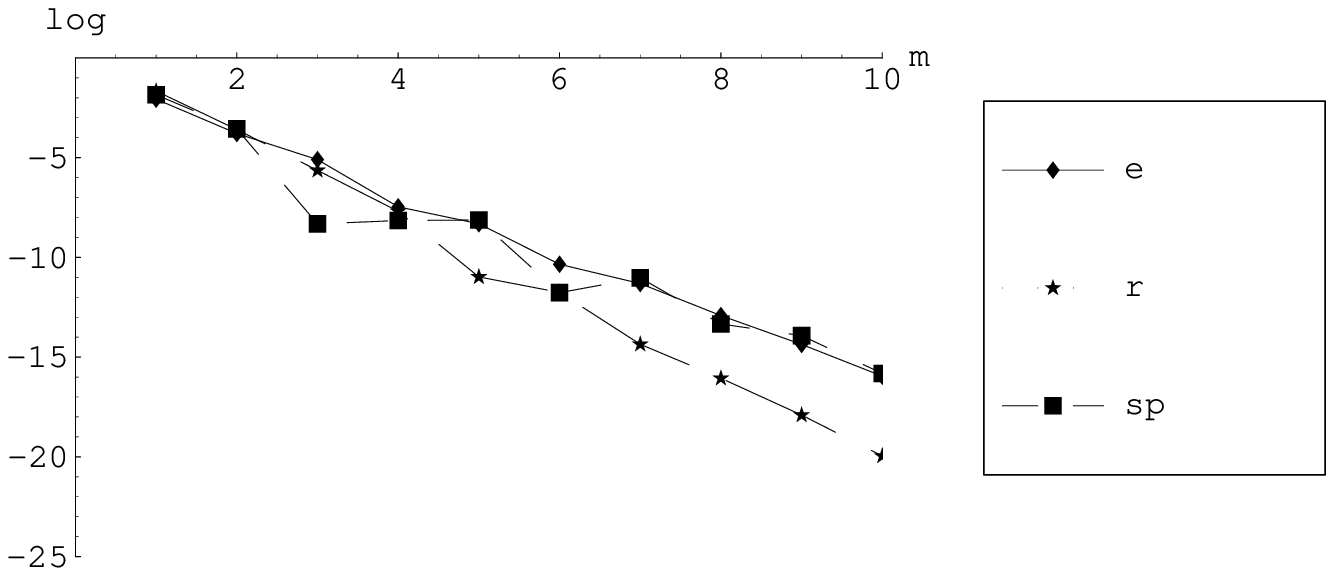}
\includegraphics{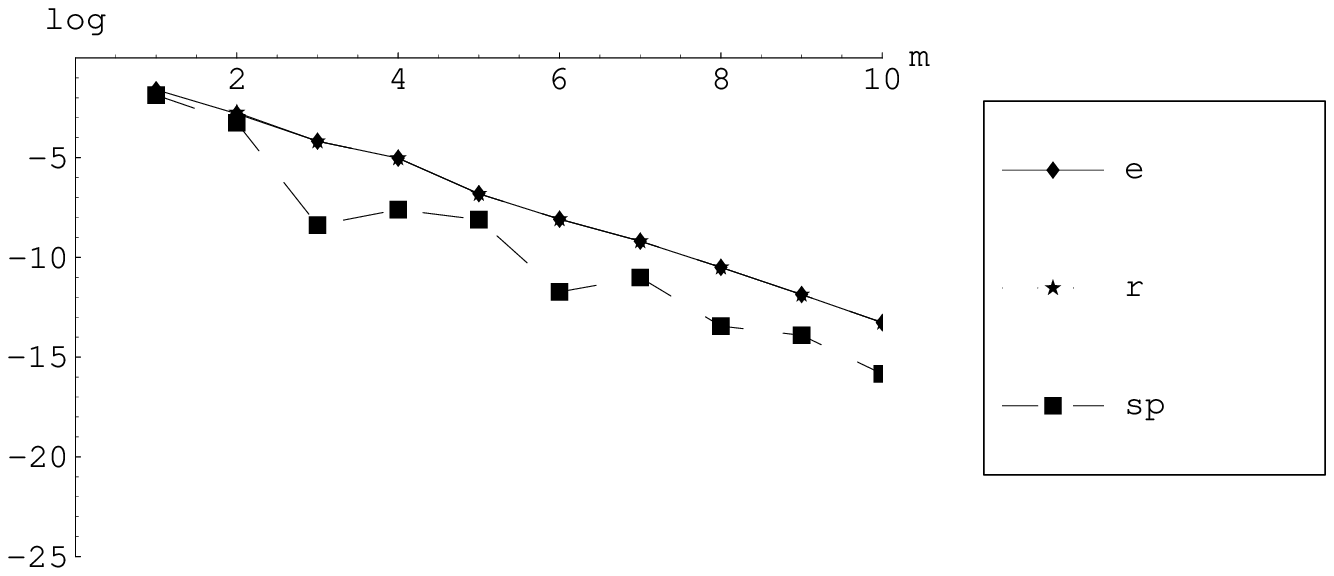}
\includegraphics{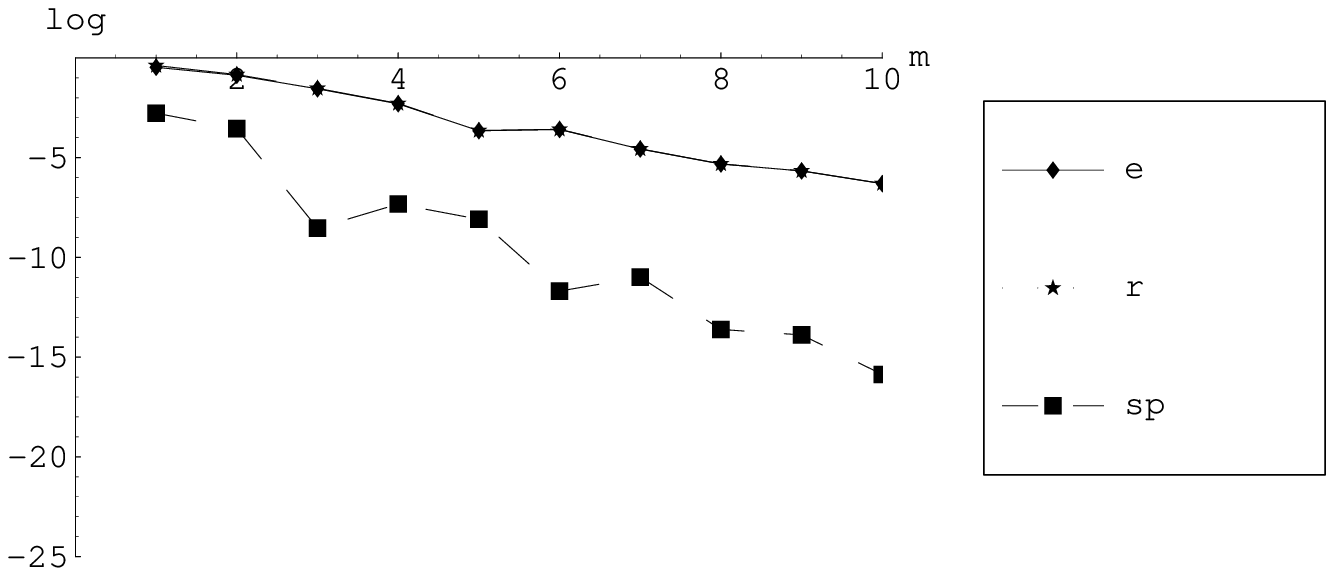}
\caption{The error \( I-T_m  \) and contribution from residues \( R_m \) and saddle points \( S_m \) for \( a = -1,0,1\)\label{graph}}
\end{figure}

We evaluate the three quantities \( I-T_m \), \( R_m \) and \( S_m \)  for \( m=1,4,...10 \), \( a=-1,0,1 \), \( t=1 \) and \( b=1 \) and presenting the results in Table \ref{err} and we plot their base \( 10 \) logarithms in Figure \ref{graph}.

We observe that for \( a \geq 0 \) the error \( I -T_m \) is determined by \( R_m \) and \(S_m \) is negligible in comparison. Furthermore, we observe that
\begin{equation}
I-T_m \sim R_m \sim \begin{cases} C \exp \left(-2m \arctan(\frac{b}{a}) \right) , & a>0 ,\\
C \exp \left(-2m \frac{\pi}{2}) \right), & a=0.
\end{cases}
\label{est2}
\end{equation}

While for \( a < 0 \) the error \( I -T_m \) is determined by \( S_m \) while \( R_m \) is negligible in comparison. Furthermore, we observe that
\begin{equation}
I-T_m \sim S_m \sim C \exp \left(-\pi m) \right).
\end{equation}

\begin{table}[tbp]
{\footnotesize \begin{tabular}{|c|c|c|c|c|c|c|c|c|c|}\hline
\multicolumn{1}{|c|} { } & \multicolumn{3}{|c|} {\(a=- 1\) } & \multicolumn{3}{|c|} {\(a=0\)} & \multicolumn{3}{|c|}{ \( a=1\)} \\ \hline
\( m \) & \(I-T_m\) & \( R_m \) & \( S_m \) & \(I-T_m\) & \( R_m \) & \( S_m \) & \(I-T_m\) & \( R_m \) & \( S_m \) \\ \hline
1&	 8.60E-3& 2.20E-2&-1.40E-2&-2.50E-2& 0.00E0&-1.40E-2&-2.50E-2& 0.00E0&-1.40E-2 \\
2&	-1.60E-4&-2.60E-4&-2.80E-4&-1.60E-3&-1.80E-3&-5.60E-4&-1.60E-3&-1.80E-3&-5.60E-4 \\
3&	 7.90E-6& 2.30E-6&-4.80E-9&-6.70E-5&-6.70E-5&-4.10E-9&-6.70E-5&-6.70E-5&-4.10E-9 \\
4&	-3.40E-8&-2.00E-8&-7.00E-9& 9.40E-6& 9.40E-6&-2.50E-8& 9.40E-6& 9.40E-6&-2.50E-8 \\
5&	-4.90E-9& 1.00E-11& 7.40E-9& 1.60E-7& 1.50E-7& 7.80E-9& 1.60E-7& 1.50E-7& 7.80E-9 \\
6&	-4.50E-11& 1.80E-12& 1.70E-12&-8.20E-9&-8.30E-9& 1.90E-12&-8.20E-9&-8.30E-9& 1.90E-12 \\
7&	 4.90E-12& 4.40E-15& 9.40E-12&-6.50E-10&-6.60E-10& 9.80E-12&-6.50E-10&-6.60E-10& 9.80E-12 \\
8&	 1.20E-13&-8.70E-17&-4.60E-14&-3.10E-11&-3.20E-11&-3.60E-14&-3.10E-11&-3.20E-11&-3.60E-14 \\
9&	-4.30E-15&-1.20E-18&-1.20E-14&-1.40E-12&-1.40E-12&-1.20E-14&-1.40E-12&-1.40E-12&-1.20E-14 \\
10&	-1.10E-16&-1.10E-20& 1.50E-16&-5.30E-14&-5.40E-14& 1.40E-16&-5.30E-14&-5.40E-14& 1.40E-16 \\ 
\hline
\end{tabular}
}
\caption{The error \( I-T_m  \) and contribution from residues \( R_m \) and saddle points \( S_m \) for \( a = -1,0,1\)\label{err}}
\end{table}

\section{Truncation error and comparison with double exponential case \label{truncerror}}
\subsection{Matching Truncation and Discretisation Errors}
We now estimate the error introduced by truncating the infinite series \( T_m \) to the finite series \( T_{n,m} \) defined by:
\begin{equation}
T_{n,m}= \frac{\pi}{m} \sum_{k=-n}^n F_h(k\frac{\pi}{m}).
\label{trunc}
\end{equation}
This error is bound by
\begin{equation}
\label{traptruncerror}
|T_m - T_{n,m}| \leq \frac{\pi}{m}  \sum_{|k|>n} |f(m \phi(k\frac{\pi}{m} ) )  \sin( m \phi(k\frac{\pi}{m} ) ) m\phi'(k\frac{\pi}{m} )|
\end{equation}
We we have that \( |f(x)| \) is bounded by some constant \( C_f \). Also, we have that
\begin{equation}
\phi'(k\frac{\pi}{m} ) = \frac{e^{k\frac{\pi}{m} }}{e^{k\frac{\pi}{m} }+1}<1.
\end{equation}
Thus, we must investigate the asymptotic behaviour of \( \sin(m \phi(k\frac{\pi}{m} )) \) as \( k \rightarrow \pm \infty \).
First, using (\ref{assminf}) we have, as \( k \rightarrow -\infty \), that
\begin{eqnarray}
\sin(m \phi(k\frac{\pi}{m} )) & = & \sin(m e^{k\frac{\pi}{m} }+O(e^{2k\frac{\pi}{m} })) \\
& = & m e^{k\frac{\pi}{m} }+O(e^{2k\frac{\pi}{m}}).
\end{eqnarray}
By (\ref{asspinf}) we have, as \( k \rightarrow \infty \), that
\begin{eqnarray}
\sin(m \phi(k\frac{\pi}{m} )) & = & \sin(k \pi +m e^{-k\frac{\pi}{m}}+O(e^{-2k\frac{\pi}{m} }))\\
& = & (-1)^k \sin(m e^{-k\frac{\pi}{m} }+O(e^{-2k\frac{\pi}{m} })) \\
& = & (-1)^k m e^{-k\frac{\pi}{m} }+O(e^{-2k\frac{\pi}{m}}).
\end{eqnarray}

Thus, from (\ref{traptruncerror}) we have that
\begin{eqnarray}
|T_m-T_{n,m}| & \leq & 2 \frac{\pi}{m}  \sum_{k=n+1}^\infty C_f(m e^{-k\frac{\pi}{m} }+O(e^{-2k\frac{\pi}{m} }))\\
& \leq & 2 \pi C_f \frac{e^{-(n+1)\frac{\pi}{m} }}{1-e^{-\frac{\pi}{m} }}\\
& \leq & 2 \pi C_f \frac{e^{-n\frac{\pi}{m} }}{\frac{\pi}{m} }.
\end{eqnarray}

Thus, if the discretisation error is, say,
\begin{equation}
I-T_m \sim C e^{-\alpha m},
\end{equation}
then we equate the exponents in the previous two equations to determine the dependence of \( m \) and \( n \), that is,
\begin{equation}
m=\sqrt{\frac{n\pi}{\alpha}}.
\end{equation}
This choice of \( m \) results in an overall error
\begin{equation}
I-T_{n,m} \sim C \sqrt{n} e^{-\sqrt{\alpha n}}.
\end{equation}

\subsection{Single exponential vs Double exponential}
\begin{figure}[tbp]
\includegraphics{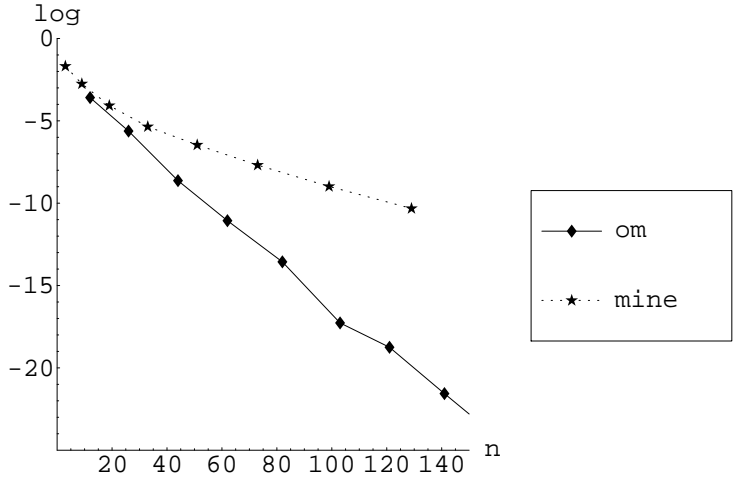}
\caption{\label{fig:entire}Our method vs Ooura and Mori's First method for \( I_1 \)}
\end{figure}

In Figure \ref{fig:entire} we present the errors in approximating the integral
\begin{equation}
I_1=\int_0^\infty \frac{\sin x}{x}\, dx
\end{equation}
using the first method of Ooura and Mori and using our method. Clearly our method results in slower convergence than theirs.

\end{document}